\pdfoutput=1
\documentclass[11pt]{amsart}
\usepackage[margin=1.15in]{geometry}
\usepackage[T1]{fontenc}
\usepackage{lmodern}
\usepackage{amssymb,amsmath,mathtools}
\usepackage{enumitem}
\usepackage{xcolor}
\usepackage{tikz}
\usetikzlibrary{arrows.meta}
\usepackage[hidelinks]{hyperref}
\hypersetup{pdftitle={On the toric Bando--Kasue--Nakajima conjecture},
  pdfauthor={Lars Andersson, Mattias Dahl, Gustav Nilsson}}

\newtheorem{theorem}{Theorem}[section]
\newtheorem{lemma}[theorem]{Lemma}
\newtheorem{proposition}[theorem]{Proposition}
\newtheorem{corollary}[theorem]{Corollary}

\newtheorem{introthm}{Theorem}

\theoremstyle{definition}
\newtheorem{definition}[theorem]{Definition}
\newtheorem{example}[theorem]{Example}
\theoremstyle{remark}
\newtheorem{remark}[theorem]{Remark}
\numberwithin{equation}{section}

\newcommand{\R}{\mathbb R}
\newcommand{\Z}{\mathbb Z}
\newcommand{\Q}{\mathbb Q}
\newcommand{\C}{\mathbb C}
\newcommand{\Tt}{T^2}                      %
\newcommand{\SU}{\mathrm{SU}}
\newcommand{\U}{\mathrm U}

\newcommand{\rb}{\bar r}                   %

\title{On the toric Bando--Kasue--Nakajima conjecture}

\author[L. Andersson]{Lars Andersson}
\address[Lars Andersson]{Department of Mathematics, KTH Royal Institute of Technology,\newline 100\,44~Stockholm, Sweden}
\email{larand@kth.se}

\author[M. Dahl]{Mattias Dahl}
\address[Mattias Dahl]{Department of Mathematics, KTH Royal Institute of Technology,\newline 100\,44~Stockholm, Sweden}
\email{dahl@kth.se}

\author[G. Nilsson]{Gustav Nilsson}
\address[Gustav Nilsson]{Mathematisches Institut, Universit\"at M\"unster,\newline Einsteinstrasse 62, 48149 M\"unster, Germany}
\email{gustav.nilsson@uni-muenster.de}

\date{November 16, 2026}

\begin{document}
\begin{abstract}
The conjecture of Bando, Kasue, and Nakajima asserts that every simply connected Ricci-flat $4$-manifold that is asymptotically locally
Euclidean (ALE), that is, asymptotic at infinity to a flat cone $\R^4/\Gamma$, is hyperk\"ahler for one of its orientations. A toric ALE instanton is a Ricci-flat ALE $4$-manifold with an effective isometric action of the $2$-torus; such a
manifold is simply connected. We prove the conjecture for the toric ALE instantons whose group at infinity $\Gamma$ is conjugate in
$\mathrm O(4)$ to a subgroup of $\SU(2)$: each of them is a multi-Eguchi--Hanson metric with collinear centres. In particular, the conjecture holds for toric
ALE instantons with $|\Gamma|\in\{1,2,3,4,6\}$. The proof uses a closed formula for toric gravitational instantons, which expresses the
metric through Hankel determinants of a positive measure on the fixed points of the torus action. With it, the problem becomes a
question about lattice points in a parallelogram determined by the end.
\end{abstract}
\maketitle

\section{Introduction}\label{sec:intro}

A gravitational instanton is a complete, non-compact, oriented, Ricci-flat Riemannian $4$-manifold whose curvature decays at least
quadratically \cite[Def.~1.1]{LiSun2026survey}. It is \emph{asymptotically locally Euclidean} (ALE) if, outside a compact set, it is
asymptotic to a flat cone $\R^4/\Gamma$, where $\Gamma\subset\mathrm O(4)$ is a finite group acting freely on the unit sphere, the
\emph{group at infinity}. The hyperk\"ahler ALE instantons are completely understood. Their group $\Gamma$ is a finite subgroup of
$\SU(2)$, and for every such $\Gamma$ Kronheimer constructed and classified them \cite{Kronheimer1989a,Kronheimer1989b}. The simplest
examples are the Eguchi--Hanson metric \cite{EguchiHanson1978}, with $\Gamma=\Z_2$, and the multi-Eguchi--Hanson metrics of Gibbons
and Hawking \cite{GibbonsHawking1978}, with $\Gamma=\Z_m$.

A conjecture of Bando, Kasue, and Nakajima \cite[Rem.~(1.8)(3)]{BKN1989} asserts that there are no others: every simply connected
Ricci-flat ALE $4$-manifold should be hyperk\"ahler for one of its two orientations. See also \cite[Conj.~11]{LiSun2026survey}, where it
is stated without simple connectivity, as the assertion that the curvature is anti-self-dual for one orientation. Gibbons had
conjectured earlier that every ALE instanton is half flat and that the group at infinity of its universal covering lies in
$\SU(2)$ \cite[Conjecture~(1)]{Gibbons1980}. A Ricci-flat metric on a simply connected manifold is half flat exactly when it is
hyperk\"ahler for one orientation, so for simply connected manifolds the two formulations agree.

Simple connectivity cannot be dropped: there are free quotients of Gibbons--Hawking metrics that are Ricci-flat K\"ahler ALE but
hyperk\"ahler for neither orientation \cite{Wright2012,Suvaina2012}.

The conjecture has been proved by Nakajima \cite[Cor.~3.4]{Nakajima1990} and the second author \cite[Thm.~6.2]{Dahl1997} when the
manifold is spin and $\Gamma\subset\SU(2)$, with the trivial spin structure on the end. Li and Sun
proved that a toric ALE instanton that is spin is hyperk\"ahler for one orientation, with no condition on $\Gamma$
\cite[Prop.~A.21]{LiSun2026toric}.

Further cases are known under additional geometric assumptions \cite{Li2023,LV2016,Ghosh2026}; without such
assumptions the conjecture is open, even for $\Gamma\subset\SU(2)$. In this paper we consider toric ALE
instantons, which carry an effective isometric action of the torus $\Tt$ (Definition~\ref{def:rods-toric}); for the literature on toric
gravitational instantons see \cite[\S\S5, 6, 7.1]{LiSun2026survey} and \cite{KL2026,AL2025}.

The group at infinity of a toric ALE instanton is cyclic, of order $m$ say, and the cross-section $S^3/\Gamma$ of its end is a lens
space. The group $\Gamma$ is conjugate in $\mathrm O(4)$ to a subgroup of $\SU(2)$ exactly when this lens space is $L(m,1)$
(Lemma~\ref{lem:rods-end}). Every toric ALE instanton is simply connected, whatever its group at infinity (Lemma~\ref{lem:basic}), so the
conjecture of Bando, Kasue, and Nakajima applies to all of them. Our main result is the following.

\begin{introthm}\label{introthm:main}
A toric ALE instanton $(M,g)$ \textup{(}Definition~\ref{def:rods-toric}\textup{)} whose group at infinity $\Gamma$ is conjugate in
$\mathrm O(4)$ to a subgroup of $\SU(2)$ is a multi-Eguchi--Hanson metric with $|\Gamma|$ collinear centres, hyperk\"ahler for one of its
orientations.
\end{introthm}

Theorem~\ref{introthm:main} assumes no spin structure, no Hermitian structure, and no topological condition on $M$. See
Theorem~\ref{thm:main} for a more detailed statement. Li and Sun proved the case $\Gamma\cong\Z_2$ \cite[Prop.~A.24]{LiSun2026toric},
as well as the spin case recalled above.

The proof rests on the solitonic form of Geroch's conjecture for toric gravitational instantons proved in \cite{GerochCF}. Its
proof gives the metric in closed form, which allows algebraic and combinatorial methods to be applied: the rod vectors of the finite
rods of a non-hyperk\"ahler instanton lie in the interior of a parallelogram in the Lie algebra of the torus that is determined by the
end \cite{N4}. What is new here is the arithmetic of the lattice points in this parallelogram. When $\Gamma$ is conjugate in
$\mathrm O(4)$ to a subgroup of $\SU(2)$, the lattice points in its interior lie on one line, too few to accommodate the rod vectors of a
non-hyperk\"ahler instanton, and the remaining hyperk\"ahler case is the multi-Eguchi--Hanson metric.

For small groups at infinity the hypothesis of Theorem~\ref{introthm:main} holds automatically. For $m\in\{1,2,3,4,6\}$ every integer
prime to $m$ is congruent to $1$ or $-1$ modulo $m$, so every lens space $L(m,q)$ is diffeomorphic to $L(m,1)$, and
Theorem~\ref{introthm:main} has the following consequence (Corollary~\ref{cor:smallgroups}).

\begin{introthm}\label{introthm:small}
A toric ALE instanton $(M,g)$ \textup{(}Definition~\ref{def:rods-toric}\textup{)} whose group at infinity $\Gamma$ has order $1$, $2$, $3$, $4$,
or $6$ is a multi-Eguchi--Hanson metric with $|\Gamma|$ collinear centres.
\end{introthm}
\section{Toric ALE instantons and rod structures}\label{sec:rods}

A toric ALE instanton is described by two kinds of data. The continuous data are the positions of the fixed points of the torus
action on the boundary of the orbit space. The discrete data are the stabilisers of the action along this boundary, a finite
sequence of lattice vectors, which we encode in a chain of integers. The group at infinity and the topology of $M$ are determined by
the chain. We recall these notions, following \cite[\S2.1]{N4} and \cite[\S2]{GerochCF}, and fix the conventions used throughout.

\subsection{Toric ALE instantons}\label{subsec:rods-toric}
Let $\Tt$ be the $2$-torus, $\mathfrak t$ its Lie algebra, and $\Lambda\subset\mathfrak t$ the lattice of the elements $X$ with
$\exp(2\pi X)=1$, so that $\Tt=\mathfrak t/2\pi\Lambda$. A \emph{$\Z$-basis} of $\Lambda$ is a basis of $\Lambda$ as a free abelian
group. The flat space $\R^4=\C^2$ carries the standard action of the torus $(S^1)^2$ by
$(z_1,z_2)\mapsto(e^{i\phi_1}z_1,e^{i\phi_2}z_2)$.

\begin{definition}\label{def:rods-toric}
A \emph{toric ALE instanton} is a complete, connected, orientable Ricci-flat Riemannian $4$-manifold $(M,g)$ with an effective
isometric action of $\Tt$, such that $(M,g)$ is ALE in the following sense. There are a compact set $K\subset M$, a finite group
$\Gamma\subset\mathrm O(4)$ that acts freely on $S^3$, and a diffeomorphism $F$ of $M\setminus K$ onto the complement of a closed ball
in $\R^4/\Gamma$ such that $|\nabla^k(g-F^*g_{\R^4})|=O(\mathfrak r^{-k-\delta})$ for some $\delta>0$ and all $k\ge0$, with the norms
and covariant derivatives of $F^*g_{\R^4}$, where $\mathfrak r$ is the Euclidean distance from the origin. The group $\Gamma$ is the
\emph{group at infinity}.
\end{definition}

\begin{remark}[invariant cross-sections]\label{rem:end-chart}
The end $E=M\setminus K$ is diffeomorphic to $(R,\infty)\times S^3/\Gamma$. The chart $F$ is not assumed to be compatible with the
action, but it can be chosen so on the universal covering $\tilde E$ of the end, to which the action lifts through a finite covering torus
$\tilde\Tt$ of $\Tt$. By the argument in the proof of \cite[Prop.~4.38]{LiSun2025}, which applies the centre-of-mass construction of
Grove and Karcher \cite{GroveKarcher1973} on the end, the lift of $F$ can be replaced by a chart $\tilde F$ that is equivariant for an
isomorphism of $\tilde\Tt$ onto $(S^1)^2$ and, after composing it with $(z_1,z_2)\mapsto(z_1,\bar z_2)$ if necessary, orientation-preserving for a given
orientation of $M$. So the end is ALE in the sense of \cite[Def.~2.9]{GerochCF} for either orientation of $M$. The chart $\tilde F$
differs from the lift of $F$, on which the deck group acts through $\Gamma\subset\mathrm O(4)$, by $O(\mathfrak r^{1-\epsilon})$ for some
$\epsilon>0$, with the corresponding bounds on derivatives, which the centre-of-mass construction provides. Averaging the function $\mathfrak r\circ\tilde F$ over the deck group gives a
$\Tt$-invariant function on $E$ that differs from $\mathfrak r\circ F$ by $O(\mathfrak r^{1-\epsilon})$, and its gradient by
$O(\mathfrak r^{-\epsilon})$. Its large level sets are \emph{invariant cross-sections} $S^3/\Gamma$ of the end, and beyond each of them the
end is a collar $[c,\infty)\times S^3/\Gamma$.
\end{remark}

\begin{lemma}\label{lem:basic}
Let $(M,g)$ be a toric ALE instanton. Then the action has finitely many fixed points, at least one; $M$ is simply connected; and all
stabilisers of the action are connected. So $(M,g)$, with either orientation, is a simply connected toric gravitational instanton
in the sense of \cite[Defs.~2.1, 2.2]{GerochCF} with an ALE end in the sense of \cite[Def.~2.9]{GerochCF}, and the results of
\cite{GerochCF,N4} apply to it.
\end{lemma}

\begin{proof}
The fixed points are isolated and lie in a compact set, since the model end contains none, so there are finitely many of them
\cite[Rem.~2.10]{GerochCF}, $n$ say. The manifold $M$ deformation retracts onto a compact $\Tt$-invariant domain $\Omega$ that contains
them and is bounded by an invariant cross-section of the end, since $M\setminus\Omega^\circ$ is a collar of $\partial\Omega$ (see
Remark~\ref{rem:end-chart}). The Euler characteristic of
a closed manifold with a torus action equals that of its fixed-point set, by \cite[Ch.~II, Thm.~5.5]{Kobayashi1972} applied to a
Killing field of an invariant metric whose flow is dense in the torus. The double of $\Omega$ carries the induced action of $\Tt$,
which is smooth for an equivariant collar of $\partial\Omega$. It has $2n$ fixed points and the Euler characteristic $2\chi(\Omega)$,
since $\partial\Omega$ is a closed $3$-manifold, so $\chi(M)=\chi(\Omega)=n$. The fundamental group of $M$
is finite, since $g$ is Ricci-flat with Euclidean volume growth \cite{LiP1986,Anderson1990}, so $b_1(M)=0$. Since $M$ is connected,
orientable, and non-compact, Poincar\'e duality gives $b_4(M)=0$ and $b_3(M)=\dim H^1_c(M;\R)$, which vanishes since $b_1(M)=0$ and the
end of $M$ is connected. Hence $n=\chi(M)=1+b_2(M)\ge1$.

When the stabilisers are connected and the action has a fixed point, simple connectivity follows from \cite[Rem.~3.9]{GerochCF}. Here
connected stabilisers are not assumed. Let $\tilde M\to M$ be the universal covering, which is finite, with the lifted metric $\tilde g$
and orientation. The Killing fields
of the action lift to complete commuting Killing fields of $\tilde g$. For $X\in\Lambda$ the time-$2\pi$ flow of the lift of $X$
covers the identity of $M$, so it is a deck transformation. It fixes the points over a fixed point of $M$, where the lifted fields
vanish, so it is the identity. Hence the lifted fields generate an effective isometric action of $\Tt$ on $\tilde M$, whose fixed
points are the finitely many points over the fixed points of $M$. Since $\tilde M$ is simply connected and non-compact, its
stabilisers are connected \cite[Lemma~2.1]{N4}. By Definition~\ref{def:rods-toric}, $|\mathrm{Rm}_g|=O(\mathfrak r^{-2-\delta})$, so $g$ has
finite curvature energy $\int_M|\mathrm{Rm}_g|^2$, and that of $\tilde g$ is $|\pi_1(M)|$ times as large. So $(\tilde M,\tilde g)$ is a
gravitational instanton in the sense of \cite[Def.~2.2]{GerochCF}; it is simply connected and toric, and its action has finitely many
fixed points, at least one. Hence the Weyl coordinates $(\tilde\rho,\tilde z)$ of a $\Z$-basis of $\Lambda$ map its orbit
space homeomorphically onto the closed half-plane $\{\tilde\rho\ge0\}$ \cite[Prop.~3.8(ii)]{GerochCF}.

A deck transformation $\gamma$ is an orientation-preserving isometry that commutes with the lifted Killing fields, hence with the
action. So it preserves the Gram matrix of the Killing fields and $\tilde\rho$, and it induces an orientation-preserving isometry of
the orbit space. Hence $\tilde z\circ\gamma$ is again a harmonic conjugate of $\tilde\rho$, so $\tilde z\circ\gamma=\tilde z+c$ for a
constant $c$. Since $\gamma$ has finite order $k$, $\tilde z=\tilde z\circ\gamma^k=\tilde z+kc$, so $c=0$. Since $(\tilde\rho,\tilde z)$ is
injective on the orbit space, $\gamma$ maps every orbit onto itself, and in particular it fixes the fixed points of the action on
$\tilde M$. A deck transformation with a fixed point is the identity, so the deck group is trivial, $M$ is simply connected, and its
stabilisers are connected \cite[Lemma~2.1]{N4}.
\end{proof}

Fix orientations of $M$ and $\mathfrak t$. For a positively oriented basis of $\mathfrak t$, let $G$ be the Gram matrix of the
corresponding Killing fields, and let $\phi^1,\phi^2$ be the dual angles. Since $g$ is Ricci-flat and the action has fixed points, the
distribution orthogonal to the orbits is integrable, and the function $\rho=(\det G)^{1/2}$ is harmonic on the interior of the orbit
space $M/\Tt$ with its quotient metric \cite[\S2.1]{GerochCF}. Where $d\rho\ne0$, a harmonic conjugate $z$ of $\rho$ gives isothermal coordinates
\cite[\S2.1]{GerochCF}. As in \cite[\S2.1]{GerochCF} and \cite[\S2.1]{N4}, we choose the harmonic conjugate for which
\begin{equation}\label{eq:rods-orient}
d\rho\wedge dz\wedge d\phi^1\wedge d\phi^2>0 .
\end{equation}
Then $(\rho,z)$ are the \emph{Weyl coordinates} of the basis, and $z$ is determined up to an additive constant. If the basis is replaced
by a basis of determinant $1$ relative to it, then $\rho$ is unchanged and $z$ changes by an additive constant. We
use $\Z$-bases of $\Lambda$, except in Section~\ref{sec:formula}. Since $M$ is simply connected (Lemma~\ref{lem:basic}), the Weyl
coordinates of a positively oriented $\Z$-basis of $\Lambda$ are global: they map the orbit space homeomorphically onto the closed
half-plane $\{\rho\ge0\}$, and its interior conformally onto the open half-plane \cite[Prop.~3.8(ii)]{GerochCF}.

\begin{figure}[ht]
\centering
\begin{tikzpicture}[x=0.85cm,y=0.85cm]
\fill[gray!8] (-0.6,0) rectangle (11.6,1.9);
\draw[->,thin] (-0.6,0) -- (11.8,0) node[right,font=\small] {$z$};
\draw[->,thin] (5.5,0) -- (5.5,2.1) node[above,font=\small] {$\rho$};
\foreach \x/\l in {1.5/w_1,3.5/w_2,6.5/w_3,9.0/w_4}{\fill (\x,0) circle (2pt); \node[below,font=\small] at (\x,-0.05) {$\l$};}
\draw[very thick,blue!70!black] (-0.6,0) -- (1.5,0);
\draw[very thick,red!70!black] (1.5,0) -- (3.5,0);
\draw[very thick,blue!70!black] (3.5,0) -- (6.5,0);
\draw[very thick,red!70!black] (6.5,0) -- (9.0,0);
\draw[very thick,blue!70!black] (9.0,0) -- (11.6,0);
\node[above,font=\small] at (0.45,0.05) {$I_0,\ v_0$};
\node[above,font=\small] at (2.5,0.05) {$I_1,\ v_1$};
\node[above,font=\small] at (4.7,0.05) {$I_2,\ v_2$};
\node[above,font=\small] at (7.75,0.05) {$I_3,\ v_3$};
\node[above,font=\small] at (10.3,0.05) {$I_4,\ v_4$};
\node[font=\footnotesize,align=center] at (8.6,1.35) {free $\Tt$-orbits};
\end{tikzpicture}
\caption{The orbit space of a toric ALE instanton with $n=4$ turning points, in Weyl coordinates. Over the rod $I_j$ the circle
generated by the rod vector $v_j$ collapses, and over a turning point the whole torus collapses to a fixed point. The preimage of the
closure of a finite rod $I_j$ is a $2$-sphere $S_j$.}\label{fig:rods-axis}
\end{figure}
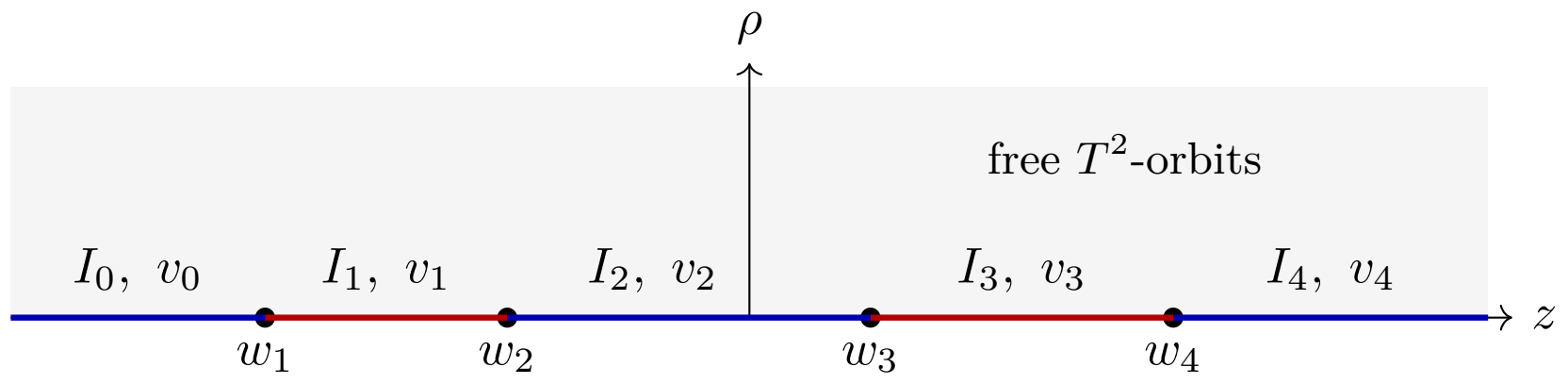

The fixed points map to $n$ points $w_1<\dots<w_n$ of the axis $\{\rho=0\}$, the \emph{turning points}. They divide the axis into the
\emph{rods} $I_0=(-\infty,w_1)$, $I_j=(w_j,w_{j+1})$ for $1\le j\le n-1$, and $I_n=(w_n,\infty)$. Over the open half-plane the
orbits are tori, over each rod the stabiliser is a circle subgroup of $\Tt$, and over a turning point it is $\Tt$ (Figure~\ref{fig:rods-axis}).
The convention \eqref{eq:rods-orient} fixes the direction of the axis, so the numbering of the rods depends on the orientations of
$M$ and $\mathfrak t$. Over the open half-plane the metric takes the Weyl--Papapetrou form
\begin{equation}\label{eq:rods-WP}
g=G_{ab}(\rho,z)\,d\phi^a\,d\phi^b+e^{2\nu(\rho,z)}\bigl(d\rho^2+dz^2\bigr),\qquad \det G=\rho^2 ,
\end{equation}
and Ricci-flatness amounts to a harmonic map equation for $\rho^{-1}G$ and a quadrature for the conformal factor $e^{2\nu}$
\cite[(2.1)--(2.4)]{GerochCF}.

\subsection{Rod vectors and chains}\label{subsec:rods-words}
Let $\det$ be the area form on $\mathfrak t$ with $\det(e,f)=1$ for every positively oriented $\Z$-basis $(e,f)$ of $\Lambda$. The
stabiliser over the rod $I_j$ is the circle subgroup generated by a primitive vector $v_j\in\Lambda$, that is, one that is not of the
form $kw$ with $w\in\Lambda$ and an integer $k\ge2$. It is the \emph{rod vector}, and it is determined up to sign. At the fixed point over $w_j$ the vectors $v_{j-1}$ and $v_j$ form a basis of $\Lambda$ \cite{OR1970}. So the
signs can be chosen, uniquely up to one common sign, such that
\begin{equation}\label{eq:rods-det}
\det(v_{j-1},v_j)=-1,\qquad 1\le j\le n .
\end{equation}
Then $(v_{j-1},v_j)$ and $(v_j,v_{j+1})$ are bases with the same determinant, so $v_{j+1}+v_{j-1}$ is an integer multiple of $v_j$:
\begin{equation}\label{eq:rods-chain}
v_{j+1}=\kappa_jv_j-v_{j-1},\qquad \kappa_j=-\det(v_{j-1},v_{j+1})\in\Z,\qquad 1\le j\le n-1 .
\end{equation}
The turning points $w_1<\dots<w_n$ together with the rod vectors $v_0,\dots,v_n$ form the \emph{rod structure} of $(M,g)$.
Since the Weyl coordinates are oriented by \eqref{eq:rods-orient}, the rod vectors and the chain of $M$ are those of
\cite{GerochCF,N4}, and the results of these papers apply to $M$ as they are stated. Finally, we identify $\Lambda$ with $\Z^2$ by the
basis $(v_1,v_0)$, which is positively oriented by \eqref{eq:rods-det}. Then
$v_0=(0,1)$, $v_1=(1,0)$, and $\det(a,b)=a_1b_2-a_2b_1$. The rod vectors are therefore determined by the integers $\kappa_j$, and we
consider such sequences on their own.

\begin{definition}\label{def:word}
A \emph{chain} is a vector $\kappa=(\kappa_1,\dots,\kappa_{n-1})\in\Z^{n-1}$, $n\ge1$. Its \emph{rod vectors} are
\begin{equation}\label{eq:rods-vectors}
v_0=(0,1),\qquad v_1=(1,0),\qquad v_{j+1}=\kappa_jv_j-v_{j-1}\quad(1\le j\le n-1),
\end{equation}
and we write $v_j=(p_j,q_j)$ and $(p,q)=(p_n,q_n)$, the last rod vector. The \emph{intersection form} of $\kappa$ is the symmetric
tridiagonal $(n-1)\times(n-1)$ matrix
\begin{equation}\label{eq:rods-Q}
Q_\kappa=\begin{pmatrix}\kappa_1&1&&\\1&\kappa_2&\ddots&\\&\ddots&\ddots&1\\&&1&\kappa_{n-1}\end{pmatrix},
\end{equation}
and the \emph{degrees} of $\kappa$ are the numbers $r$ and $\rb$ of negative and of positive eigenvalues of $Q_\kappa$. We call $r$
\emph{the degree} of $\kappa$.\footnote{Li and
Sun define the degree of a rod structure as a winding number of its rod lines \cite[(3.2), Def.~3.11]{LiSun2025}. For an ALE end their
degree of $(M,g)$ equals $\rb=b_2^+(M)$ (Lemma~\ref{lem:rods-degree}, Remark~\ref{rem:LiSun-degree}). Our $r$ is the order of the
Hankel determinants in the closed formula, the degree index of \cite{GerochCF,N4}; this was proved in \cite[Thm.~2.4(i)]{N4} (see
Lemma~\ref{lem:sol-degree} for a direct argument).} The
\emph{chain of a toric ALE instanton} $(M,g)$, with orientations of $M$ and $\mathfrak t$, is the chain defined by
\eqref{eq:rods-chain}, and its degrees are the degrees of $(M,g)$.
\end{definition}

By induction, the rod vectors of a chain satisfy \eqref{eq:rods-det}, and $p_j=-\det(v_0,v_j)$. In particular
$p=-\det(v_0,v_n)$. Our $p$ is minus the $p$ of \cite{N4}, and from \cite{N4} we use only statements about $|p|$. The chain determines $M$ with its torus action up to equivariant diffeomorphism \cite{OR1970},
\cite[proof of Thm.~1]{Nilsson2024}, and it has a direct topological meaning \cite[\S2.1]{N4}. The preimage of the closure of a finite
rod $I_j$ is an embedded $2$-sphere $S_j\subset M$ (Figure~\ref{fig:rods-axis}). Consecutive spheres meet transversally in one fixed
point, spheres that are not consecutive are disjoint, the classes of $S_1,\dots,S_{n-1}$ form a basis of $H_2(M;\Z)$, and by \cite[(2.6)]{N4}
\begin{equation}\label{eq:rods-selfint}
S_j\cdot S_j=\kappa_j,\qquad 1\le j\le n-1 .
\end{equation}
So $Q_\kappa$ is the intersection form of $M$ in the basis $(S_j)$, for suitable orientations of the spheres, and
the degrees are $r=b_2^-(M)$ and $\rb=b_2^+(M)$. Reversing the orientation of $M$ replaces $\kappa$ by $-(\kappa_{n-1},\dots,\kappa_1)$
\cite[\S2.1]{N4} and exchanges the two degrees. Reversing the orientation of $\mathfrak t$ replaces $\kappa$ by
$(\kappa_{n-1},\dots,\kappa_1)$ and leaves the degrees unchanged.

\begin{remark}[orientation conventions]\label{rem:rods-HJ}
With the opposite convention $d\rho\wedge dz\wedge d\phi^1\wedge d\phi^2<0$, the chain assigned to $M$ is $-(\kappa_{n-1},\dots,\kappa_1)$,
so its entries are minus the self-intersection numbers of the spheres, in reverse order. In that convention the chain of the minimal resolution of a cyclic
quotient singularity $\C^2/\Gamma_{m,a}$ (Section~\ref{subsec:rods-end}), with its complex orientation, consists of the entries of the
Hirzebruch--Jung continued fraction of $m/a$, all of which are $\ge2$. This is the usual convention in algebraic geometry
\cite{Hirzebruch1953}, \cite[\S2.6]{Fulton1993}, and the
multi-Eguchi--Hanson metrics have the chain $(2,\dots,2)$ in their hyperk\"ahler orientation. We follow the convention of
\cite{GerochCF,N4}, in which they have the chain $(-2,\dots,-2)$ (Example~\ref{ex:mEH}).
\end{remark}

\subsection{The end}\label{subsec:rods-end}
The two semi-infinite rods determine the end of $M$. We recall the arithmetic needed to describe it. For an integer $m\ge1$, two
integers $a$ and $a'$ are \emph{congruent modulo $m$}, $a\equiv a'\pmod m$, if $m$ divides $a-a'$. A \emph{unit modulo $m$} is an
integer prime to $m$. The units are exactly the integers $a$ for which $aa'\equiv1\pmod m$ for some integer $a'$, an inverse of $a$
modulo $m$, and the number of units among $1,\dots,m$ is Euler's function $\varphi(m)$. The integers $1$ and $-1$ are always units,
and they are the only ones up to congruence exactly when $\varphi(m)\le2$, that is, when $m\in\{1,2,3,4,6\}$. For $m=5$, for
instance, the unit $2$ is congruent to neither $1$ nor $-1$.

For a unit $a$ modulo $m$ let $\Gamma_{m,a}\subset\U(2)$ be the cyclic group generated by $\operatorname{diag}(\zeta,\zeta^a)$,
$\zeta=e^{2\pi i/m}$. It acts freely on $S^3$, and the quotient is the lens space $L(m,a)$. The lens spaces $L(m,a)$ and $L(m,a')$ are
homeomorphic if and only if $a'\equiv\pm a$ or $aa'\equiv\pm1\pmod m$ \cite{Brody1960}, and then they are diffeomorphic. Since
$\det(v_{n-1},v_n)=-1$, the coordinates $p$ and $q$ of the last rod vector are coprime, so that $q$ is a unit modulo $|p|$ if $p\ne0$.

\begin{lemma}\label{lem:rods-end}
Let $(M,g)$ be a toric ALE instanton with the chain $\kappa$, and let $\Gamma\subset\mathrm O(4)$ be its group at infinity. Then
$p\ne0$, $\Gamma\cong\Z_{|p|}$, the cross-section $S^3/\Gamma$ of the end is diffeomorphic to $L(|p|,q)$, and $\Gamma$ is conjugate in
$\mathrm O(4)$ to a subgroup of $\SU(2)$ if and only if $q\equiv\pm1\pmod{|p|}$.
\end{lemma}

\begin{proof}
Let $\Sigma\cong S^3/\Gamma$ be an invariant cross-section of the end (see Remark~\ref{rem:end-chart}). Its orbit
space is a compact $1$-manifold in the closed half-plane \cite[\S2]{OR1970}. It is not a circle, since then $\Sigma$ would be a torus bundle
over the circle, whose fundamental group is infinite. So it is an arc with its end points on the axis. Since $\Sigma$ separates the fixed
points from infinity, the arc separates the turning points from infinity in the half-plane, so its end points lie on the rods $I_0$ and
$I_n$, with the stabilisers generated by $v_0$ and $v_n$, and the stabiliser is trivial in between. If $p=0$, then $v_n=\pm v_0$, and $\Sigma$ is diffeomorphic to $S^1\times S^2$, whose fundamental group is infinite. Since
$\Gamma$ is finite, $p\ne0$. A $3$-manifold with such an action is determined up to equivariant diffeomorphism by the lattice and the
two stabilisers \cite[\S2]{OR1970}, and comparison with a linear action on a lens space shows that $\Sigma$ is $L(|p|,q)$
\cite[Lemma~3]{Nilsson2024}. Hence $\Gamma=\pi_1(\Sigma)\cong\Z_{|p|}$.

No non-trivial element of $\Gamma$ has the eigenvalue $1$, since $\Gamma$ acts freely on $S^3$. Every element of $\mathrm O(4)\setminus\mathrm{SO}(4)$
has the eigenvalue $1$, so $\Gamma\subset\mathrm{SO}(4)$, and $\Gamma$ is conjugate in $\mathrm O(4)$ to $\Gamma_{|p|,a}$ for a unit $a$. By
the classification of lens spaces, $a\equiv\pm q$ or $aq\equiv\pm1$. A generator of $\Gamma_{|p|,a}$ rotates two orthogonal planes by
the angles $2\pi/|p|$ and $2\pi a/|p|$, while an element of $\SU(2)$ rotates two orthogonal planes by opposite angles. So
$\Gamma_{|p|,a}$ is conjugate in $\mathrm O(4)$ into $\SU(2)$ if and only if $a\equiv\pm1$, that is, if and only if $q\equiv\pm1$.
\end{proof}

We say that $M$ has the \emph{end} $L(m,a)$ if the cross-section $S^3/\Gamma$ is diffeomorphic to $L(m,a)$, without reference to an
orientation. For a chain with $p\ne0$ we likewise call $L(|p|,q)$ its \emph{end}. By Lemma~\ref{lem:rods-end} and the classification of
lens spaces, the following three conditions are therefore
equivalent: $\Gamma$ is conjugate in $\mathrm O(4)$ to a subgroup of $\SU(2)$, the end is $L(m,1)$ with $m=|\Gamma|$, and $q\equiv\pm1\pmod{|p|}$.

\begin{example}\label{ex:mEH}
The \emph{multi-Eguchi--Hanson metrics} with $n$ collinear centres are the hyperk\"ahler ALE metrics obtained from the
Gibbons--Hawking ansatz with $n$ centres on a line in $\R^3$ \cite{GibbonsHawking1978}, \cite[Ex.~4.3]{LiSun2026survey},
\cite[(2.10)]{N4}. For $n=1$ the metric is flat, and for $n=2$ it is the Eguchi--Hanson metric \cite{EguchiHanson1978}. They are toric
ALE instantons whose turning points are the centres \cite[\S11.2]{GerochCF}. In their hyperk\"ahler orientation the chain is
$(-2,\dots,-2)$, and $Q_\kappa$ is minus the Cartan matrix of $A_{n-1}$, so that $r=n-1$ and $\rb=0$ \cite[Rem.~3.6]{N4}. With the
opposite orientation the chain is $(2,\dots,2)$, and the two degrees are exchanged. The rod vectors and the end $L(n,1)$ are given in
Table~\ref{tab:chains}.
\end{example}

Table~\ref{tab:chains} collects the data of the chains that serve as examples in this paper: the chains of Example~\ref{ex:mEH}, the
chain $(2,2,2)$ of the multi-Eguchi--Hanson metrics with four centres, with the orientation opposite to the hyperk\"ahler one, and the
chains $(1,-1,1,-1)$ and $(2,-2)$, whose end $L(5,2)$ is not covered by our argument (Example~\ref{ex:diamond}). The rod vectors follow
from \eqref{eq:rods-vectors}, the ends from Lemma~\ref{lem:rods-end} and the classification of lens spaces, and the degrees from the
eigenvalues of $Q_\kappa$, or more quickly from Lemma~\ref{lem:rods-degree} below.

\begin{table}[ht]
\centering
\small
\renewcommand{\arraystretch}{1.2}\setlength{\tabcolsep}{4pt}
\begin{tabular}{@{}lclccc@{}}
\hline
chain $\kappa$ & $n$ & rod vectors $v_0,\dots,v_n$ & $(p,q)$ & end & $(r,\rb)$\\
\hline
$(2,\dots,2)$ & $n$ & $v_j=(j,1-j)$ & $(n,1-n)$ & $L(n,1)$ & $(0,n-1)$\\
$(-2,\dots,-2)$ & $n$ & $v_j=(-1)^{j+1}(j,j-1)$ & $(-1)^{n+1}(n,n-1)$ & $L(n,1)$ & $(n-1,0)$\\
$(2,2,2)$ & $4$ & $(0,1)$, $(1,0)$, $(2,-1)$, $(3,-2)$, $(4,-3)$ & $(4,-3)$ & $L(4,1)$ & $(0,3)$\\
$(1,-1,1,-1)$ & $5$ & $(0,1)$, $(1,0)$, $(1,-1)$, $(-2,1)$, $(-3,2)$, $(5,-3)$ & $(5,-3)$ & $L(5,2)$ & $(2,2)$\\
$(2,-2)$ & $3$ & $(0,1)$, $(1,0)$, $(2,-1)$, $(-5,2)$ & $(-5,2)$ & $L(5,2)$ & $(1,1)$\\
\hline
\end{tabular}
\caption{The chains used as examples, with $n$ the number of turning points, the rod vectors \eqref{eq:rods-vectors}, the last rod
vector $(p,q)$, the end (Lemma~\ref{lem:rods-end}), and the degrees $r$ and $\rb$ (Definition~\ref{def:word}). In the first two rows
$0\le j\le n$.}\label{tab:chains}
\end{table}
\section{Chains and their degrees}\label{sec:degree}

The degree $r$ of a chain was defined in Definition~\ref{def:word} as the number of negative eigenvalues of the intersection form
$Q_\kappa$ of \eqref{eq:rods-Q}. It can also be read off from the rotation of the rod vectors, and this description identifies $\rb$
with the degree of Li and Sun. By \eqref{eq:rods-det}, each $v_j$ is
obtained from $v_{j-1}$ by a clockwise rotation through an angle $\theta_j\in(0,\pi)$, followed by a dilation. Here clockwise refers to
the orientation of $\Lambda=\Z^2$, and the angles are measured in the Euclidean structure of $\R^2$. The \emph{total turning} of the
rod vectors is
\begin{equation}\label{eq:degree-Theta}
\Theta=\theta_1+\dots+\theta_n ,
\end{equation}
the clockwise angle from $v_0$ to $v_n$ swept out by the rod vectors. We also put $\varepsilon=(-1)^r$.

\begin{lemma}\label{lem:rods-degree}
Let $\kappa$ be a chain with $p\ne0$, with the intersection form $Q_\kappa$ of \eqref{eq:rods-Q}, the degrees $r$ and $\rb$ of
Definition~\ref{def:word}, and the total turning $\Theta$ of \eqref{eq:degree-Theta}. Then $\det Q_\kappa=p$, so that $Q_\kappa$ is
non-degenerate and $r+\rb=n-1$. Moreover $\Theta\notin\pi\Z$, and
\begin{equation}\label{eq:rods-Theta}
r=\lfloor\Theta/\pi\rfloor
\end{equation}
is the number of multiples of $\pi$ in the interval $(0,\Theta)$. Equivalently, $r$ is the number of sign changes in the sequence
$(p_1,\dots,p_n)$ of first coordinates of the rod vectors, where a vanishing term counts as one sign change between its neighbours.
In particular $p$ has the sign $\varepsilon=(-1)^r$.
\end{lemma}

\begin{proof}
Let $D_j$ be the leading principal $j\times j$ minor of $Q_\kappa$, with $D_0=1$. Expanding along the last row gives $D_1=\kappa_1$ and
$D_j=\kappa_jD_{j-1}-D_{j-2}$, while \eqref{eq:rods-vectors} gives $p_0=0$, $p_1=1$, and $p_{j+1}=\kappa_jp_j-p_{j-1}$. By induction,
\begin{equation}\label{eq:rods-minors}
D_j=p_{j+1},\qquad 0\le j\le n-1 .
\end{equation}
In particular $\det Q_\kappa=D_{n-1}=p\ne0$.

Let $\Theta_j=\theta_1+\dots+\theta_j$ be the clockwise angle from $v_0$ to $v_j$, so that $\Theta_n=\Theta$. Since $v_0=(0,1)$,
$p_j=|v_j|\sin\Theta_j$. The angles $\Theta_j$ increase in steps less than $\pi$. So a vanishing $p_j$ with $0<j<n$ has neighbours $p_{j-1},p_{j+1}$ of opposite signs, and
two consecutive terms do not vanish. Since $\Theta_1\in(0,\pi)$ and $p_n=p\ne0$, the number $r'$ of sign changes in
$(p_1,\dots,p_n)$, counted as in the statement, is the number of multiples of $\pi$ in $(0,\Theta)$, and
$\Theta\notin\pi\Z$.

It remains to show that $r=r'$. We use Jacobi's rule for the signature of a real symmetric matrix, with Gundelfinger's
supplement for a vanishing minor \cite[Ch.~X, \S3, Thm.~2 and Note~1]{Gantmacher1959}. For a symmetric tridiagonal matrix with non-zero entries next to
the diagonal and non-zero determinant, the leading principal minors $(D_0,\dots,D_{n-1})$ have no two consecutive zeros, and a
vanishing minor has neighbours of opposite signs, since $D_{j+1}=-D_{j-1}$ when $D_j=0$. The rules then say that the number of
negative eigenvalues is the number of sign changes in this sequence, where a vanishing minor counts as one sign change. By
\eqref{eq:rods-minors} this sequence is $(p_1,\dots,p_n)$, with the same convention for vanishing terms, so $r=r'$. Finally
$p=|v_n|\sin\Theta$ and $\Theta\in(r\pi,(r+1)\pi)$, so $p$ has the sign $(-1)^r$.
\end{proof}

For the chains of Table~\ref{tab:chains} the lemma can be checked directly. The rod vectors of the chain $(1,-1,1,-1)$ have the first
coordinates $1,1,-2,-3,5$, which change sign twice, so $r=2$. The chains $(2,\dots,2)$ have the positive first coordinates $1,2,\dots,n$, so $r=0$, while
the chains $(-2,\dots,-2)$ have the first coordinates $(-1)^{j+1}j$, so $r=n-1$.

\begin{remark}[chains as walks in the Farey graph]\label{rem:farey}
The vertices of the Farey graph are the points of $\Q\cup\{\infty\}$, and two fractions $a_2/a_1$ and $b_2/b_1$ in lowest terms are
adjacent if $|a_1b_2-a_2b_1|=1$ \cite[Ch.~1--3]{Hatcher2022}. By \eqref{eq:rods-det}, consecutive rod vectors form a basis of $\Z^2$.
So their slopes $x_j=q_j/p_j$, with $\frac10=\infty$, are fractions in lowest terms, and consecutive slopes are adjacent. Hence a chain
is a walk $x_0=\infty$, $x_1=0$, $x_2,\dots,x_n=q/p$ in the Farey graph, that is, a sequence of vertices in which consecutive ones are
adjacent, from the slope of $v_0$ to that of $v_n$. The slope decreases along each step on the circle $\R\cup\{\infty\}$, since
$v_{j+1}$ is obtained from $v_j$ by a clockwise rotation through an angle in $(0,\pi)$.
By Lemma~\ref{lem:rods-degree}, $r$ is the number of steps from $x_j$ to $x_{j+1}$, $1\le j\le n-1$, with $p_jp_{j+1}<0$ or
$p_{j+1}=0$, that is, of the steps after the first that cross or reach $\infty$. For the chain $(1,-1,1,-1)$ the walk is
$\infty,0,-1,-\frac12,-\frac23,-\frac35$, and these are the steps from $-1$ to $-\frac12$ and from $-\frac23$ to $-\frac35$.
\end{remark}

\begin{remark}[the degree of Li and Sun]\label{rem:LiSun-degree}
For a toric ALE instanton, $\rb$ is the degree that Li and Sun assign to the rod structure of $(M,g)$ \cite[(3.2), Def.~3.11]{LiSun2025},
and $r$ is the degree of the rod structure of $M$ with the opposite orientation. We use the description of their degree in
\cite[\S2.1]{N4}, in the orientation convention \eqref{eq:rods-orient}. Each of the lines $\R v_0,\dots,\R v_n$ spanned by the rod vectors is rotated
counter-clockwise through an angle in $(0,\pi)$ onto the next, and finally $\R v_n$ onto $\R v_0$. Since $v_j$ is obtained from $v_{j-1}$ by
a clockwise rotation through $\theta_j$, the first $n$ angles are $\pi-\theta_1,\dots,\pi-\theta_n$. The last angle, counter-clockwise
from $\R v_n$ to $\R v_0$, is the clockwise angle from $\R v_0$ to $\R v_n$, which is $\Theta-r\pi\in(0,\pi)$ by Lemma~\ref{lem:rods-degree}, as $p\ne0$
by Lemma~\ref{lem:rods-end}.
The sum is $(n-r)\pi$,
and for an ALE end the degree is this sum divided by $\pi$, minus $1$, that is, $n-1-r=\rb$. Reversing the orientation of $M$ exchanges
$r$ and $\rb$ (Section~\ref{subsec:rods-words}), which gives the second statement.
\end{remark}
\section{The closed formula and the diamond inequality}\label{sec:formula}

The proof of the main theorem needs one consequence of the analytic theory of toric gravitational instantons: for a toric ALE instanton whose
degree $r$ satisfies $1\le r\le n-2$, the interior rod vectors $v_1,\dots,v_{n-1}$ lie in the open parallelogram with the vertices
$\pm v_0$ and $\pm v_n$, the diamond of the end. This inequality is due to \cite{N4}. It follows from the closed formula for the
metric \cite{GerochCF}, which expresses the rod vectors through Hankel determinants of a positive measure on the turning points, and
from Heine's formula, which writes these determinants as sums of positive terms. We recall Hankel determinants, state the closed
formula, and include the short proof of the inequality for completeness.

\subsection{Hankel determinants and Heine's formula}\label{subsec:sol-hankel}
Let $\lambda=\sum_{k=1}^N\lambda_k\delta_{x_k}$ be a measure with complex weights $\lambda_k$ on distinct real points $x_1,\dots,x_N$. Its
moments are $h_i(\lambda)=\sum_k\lambda_kx_k^i$, $i\ge0$, and its \emph{Hankel determinants} are
\begin{equation}\label{eq:heine-hankel}
H_s(\lambda)=\det\bigl(h_{i+l}(\lambda)\bigr)_{i,l=0}^{s-1}\quad(s\ge1),\qquad H_0(\lambda)=1 .
\end{equation}
For a set $S\subset\{1,\dots,N\}$ the \emph{Heine term} is
\begin{equation}\label{eq:heine-term}
P_S(\lambda)=\prod_{k\in S}\lambda_k\prod_{\substack{k,l\in S\\k<l}}(x_k-x_l)^2 ,
\end{equation}
with empty products equal to $1$.

\begin{lemma}[Heine's formula]\label{lem:heine}
For every $s\ge0$,
\begin{equation}\label{eq:heine}
H_s(\lambda)=\sum_{\substack{S\subset\{1,\dots,N\}\\|S|=s}}P_S(\lambda).
\end{equation}
\end{lemma}

\begin{proof}
For $s=0$ both sides equal $1$. For $s\ge1$ let $W$ be the $N\times s$ matrix with the entries $W_{kl}=x_k^l$, $0\le l\le s-1$. Then
$W^T\operatorname{diag}(\lambda_1,\dots,\lambda_N)W=(h_{i+l}(\lambda))_{i,l}$, and the Cauchy--Binet formula gives
$H_s(\lambda)=\sum_{|S|=s}\prod_{k\in S}\lambda_k\det(W_{S})^2$, where $W_S$ is the $s\times s$ submatrix of the rows in $S$. This is a
Vandermonde matrix, with $\det(W_S)^2=\prod_{k<l\in S}(x_k-x_l)^2$.
\end{proof}

Formula \eqref{eq:heine} goes back to Heine \cite{Heine1878}. See also \cite[Ch.~II]{Szego1975} and \cite[Lemma~5.1]{GerochCF}. In
particular, if all weights are positive, then $H_s(\lambda)>0$ for $0\le s\le N$, and $H_s(\lambda)=0$ for $s>N$.

Now let $\mu=\sum_{k=1}^n\omega_k\delta_{w_k}$ be a positive measure on $n$ points $w_1<\dots<w_n$, so that $\omega_k>0$. For $0\le j\le n$ let
\begin{equation}\label{eq:sol-muj}
\mu_j=-\sum_{k\le j}\omega_k\delta_{w_k}+\sum_{k>j}\omega_k\delta_{w_k}
\end{equation}
be $\mu$ with the signs of its first $j$ weights reversed. By Heine's formula, $H_s(\mu_j)$ is a sum of the positive Heine terms $P_S(\mu)$
with signs, and this bounds its size. For two points, for example, $H_1(\mu_1)=\omega_2-\omega_1$ lies strictly between $-H_1(\mu)$ and
$H_1(\mu)$.

\begin{lemma}[{\cite[Lemma~3.1]{N4}}]\label{lem:sol-signs}
Let $\mu$ be a positive measure on $n\ge1$ points and $0\le s\le n$. Then $H_s(\mu)>0$, and the following hold.
\begin{enumerate}[leftmargin=*,label=\textup{(\roman*)}]
\item $H_s(\mu_n)=(-1)^sH_s(\mu)$, and $H_n(\mu_j)=(-1)^jH_n(\mu)$ for $0\le j\le n$.
\item If $1\le s\le n-1$, then $|H_s(\mu_j)|<H_s(\mu)$ for $0<j<n$.
\end{enumerate}
\end{lemma}

\begin{proof}
By \eqref{eq:heine-term}, $P_S(\mu_j)=(-1)^{|S\cap\{1,\dots,j\}|}P_S(\mu)$, and $P_S(\mu)>0$. So by \eqref{eq:heine}, $H_s(\mu_j)$ is the sum of
the positive terms $P_S(\mu)$, $|S|=s$, with these signs. For $j=n$ every set of $s$ atoms has the sign $(-1)^s$, and for $s=n$ the only
set is $\{1,\dots,n\}$, with the sign $(-1)^j$. This proves (i). The triangle inequality gives $|H_s(\mu_j)|\le H_s(\mu)$, with
equality only if all signs agree. Let $0<j<n$ and $1\le s\le n-1$, and choose a set $T\subset\{2,\dots,n-1\}$ of $s-1$ atoms. The sets
$T\cup\{1\}$ and $T\cup\{n\}$ carry opposite signs, since the first atom is reversed in $\mu_j$ and the last one is not. Hence the
inequality is strict, which proves (ii).
\end{proof}

\subsection{The closed formula}\label{subsec:sol-rep}
The closed formula is built from the following functions. Let $w_1<\dots<w_n$ be points of the axis, $\mu=\sum_k\omega_k\delta_{w_k}$ a
positive measure, and $0\le s\le n-1$. For $\rho\ge0$ and $(\rho,z)$ not one of the points $(0,w_k)$ put
\begin{equation}\label{eq:B-rvarpi}
R_k=\bigl(\rho^2+(z-w_k)^2\bigr)^{1/2},\qquad \varpi_k=\frac{w_k-z-i\rho}{R_k},
\end{equation}
so that $|\varpi_k|=1$. Let $\varpi\mu$, $\mathbf R\mu$, and $\mu/\mathbf R$ be the measures on the points $w_k$ with the weights
$\varpi_k\omega_k$, $R_k\omega_k$, and $\omega_k/R_k$, and let
\begin{equation}\label{eq:B-XYF}
\mathcal X=H_s(\varpi\mu),\qquad \mathcal Y=H_{s+1}(\varpi\mu),\qquad F_1=H_s(\mathbf R\mu),\qquad F_2=H_{s+1}(\mu/\mathbf R).
\end{equation}
By Heine's formula, $F_1>0$ and $F_2>0$, since $s+1\le n$. On the rod $I_j$, $0\le j\le n$, we have $\rho=0$, $\varpi_k=-1$ for $k\le j$,
and $\varpi_k=1$ for $k>j$, so that $\varpi\mu=\mu_j$ there. For $0\le j\le n$ put
\begin{equation}\label{eq:B-Vj}
V_j=\bigl(H_s(\mu_j),\,H_{s+1}(\mu_j)\bigr)\in\R^2 .
\end{equation}

\begin{theorem}[{\cite[Thm.~10.1, (10.1)]{GerochCF}, in the form of \cite[Thm.~2.4]{N4}}]\label{thm:soliton}
Let $(M,g)$ be a toric ALE instanton with $n$ turning points and the rod vectors $v_0,\dots,v_n$ of
Section~\ref{subsec:rods-words}. Then there are an integer $s$ with $0\le s\le n-1$, a positively oriented basis $(f_0,f_1)$ of
$\mathfrak t$, a positive measure $\mu=\sum_k\omega_k\delta_{w_k}$ on the turning points $w_1<\dots<w_n$, in the Weyl coordinates
$(\rho,z)$ of this basis with $z$ chosen as in \eqref{eq:rods-orient}, and a constant $\mathfrak c>0$ with the following properties,
where $\mathcal X$, $\mathcal Y$, $F_1$, and $F_2$ are defined by \eqref{eq:B-rvarpi} and \eqref{eq:B-XYF}, and $V_j$ by
\eqref{eq:sol-muj} and \eqref{eq:B-Vj}, with the Hankel determinants $H_s$ and $H_{s+1}$ of the orders $s$ and $s+1$.
\begin{enumerate}[leftmargin=*,label=\textup{(\roman*)}]
\item With the angles $\psi^0,\psi^1$ dual to $(f_0,f_1)$,
\begin{equation}\label{eq:sol-closed}
g=\frac{\bigl|\mathcal Y\,d\psi^0-\mathcal X\,d\psi^1\bigr|^2}{F_1F_2}+\mathfrak c\,F_1F_2\,(d\rho^2+dz^2).
\end{equation}
\item After a change of the common sign of the rod vectors, which does not change the chain, the components of $v_j$ with respect to
$(f_0,f_1)$ are
\begin{equation}\label{eq:sol-rodvectors}
v_j=\mathfrak c^{1/2}\,V_j,\qquad 0\le j\le n .
\end{equation}
\end{enumerate}
\end{theorem}

Theorem~\ref{thm:soliton} is the detailed form of the theorem that toric gravitational instantons are gravitational solitons
\cite[Thm.~1]{GerochCF}, for an ALE end (for a Pfaffian form of the Belinski--Zakharov solitons see also
\cite[\S\S2.3--2.4]{Aksteiner2026}). Its hypotheses hold by Lemma~\ref{lem:basic}. The soliton data of \cite{GerochCF,N4} contain a
further number $\gamma\ge0$, the weight of an atom at infinity, which vanishes exactly for ALE ends \cite[\S2.3]{N4}, so that
\cite[(2.16)]{N4} reduces to \eqref{eq:B-XYF} and \eqref{eq:sol-closed}. Items (i) and (ii) are
items (ii) and (iii) of \cite[Thm.~2.4]{N4}, and $s$, the order of the Hankel determinants, is the degree index of \cite{GerochCF,N4}.
Since $(f_0,f_1)$ is positively oriented, its Weyl coordinates differ from those of a positively oriented $\Z$-basis of
$\Lambda$ by a positive factor and a translation, so the turning points, the rods, and their numbering are those of
Section~\ref{subsec:rods-toric}. Item (ii) says that along each rod the generator of the stabiliser is given by Hankel determinants of
the measure $\mu$ with some signs reversed. Item (i) of \cite[Thm.~2.4]{N4} relates the degree index to the degree of Li and Sun, and with
Remark~\ref{rem:LiSun-degree} it identifies the degree index with our degree $r$. The following lemma gives a direct argument from
item (ii).

\begin{lemma}\label{lem:sol-degree}
The integer $s$ of Theorem~\ref{thm:soliton} equals the degree $r$ of $(M,g)$ of Definition~\ref{def:word}.
\end{lemma}

\begin{proof}
Let
\begin{equation}\label{eq:sol-normfamily}
\tilde V_j=\Bigl(\frac{H_s(\mu_j)}{H_s(\mu)},\ \frac{H_{s+1}(\mu_j)}{H_{s+1}(\mu)}\Bigr),\qquad 0\le j\le n ,
\end{equation}
be the normalised Hankel family of \cite[\S5]{GerochCF}. It is the image of $(V_j)$ under a diagonal linear map with positive
entries. By \eqref{eq:sol-rodvectors}, and since $(f_0,f_1)$ is positively oriented, $\det(v_i,v_j)$ is therefore a positive multiple of
the determinant $\det(\tilde V_i,\tilde V_j)$ of the components, for $0\le i,j\le n$. In particular $\det(\tilde V_{j-1},\tilde V_j)<0$ by
\eqref{eq:rods-det}, so $\tilde V_j$ is obtained from $\tilde V_{j-1}$ by a clockwise rotation through an angle in $(0,\pi)$, followed by a
dilation, and these $n$ angles add up to $(s+\frac12)\pi$ \cite[Prop.~5.5]{GerochCF}. By Lemma~\ref{lem:sol-signs}(i),
$\tilde V_0=(1,1)$ and $\tilde V_n=(-1)^s(1,-1)$, so $\det(\tilde V_0,\tilde V_n)\ne0$. The argument of the proof of
Lemma~\ref{lem:rods-degree}, with $\tilde V_0$ in place of $v_0$, shows that the number of sign changes in the sequence
$(\det(\tilde V_0,\tilde V_1),\dots,\det(\tilde V_0,\tilde V_n))$, counted as there, is the number of multiples of $\pi$ in
$(0,(s+\frac12)\pi)$, that is, $s$. This sequence has the signs of $(-p_1,\dots,-p_n)$, so by Lemma~\ref{lem:rods-degree} the number of
its sign changes is $r$. Hence $s=r$.
\end{proof}

\subsection{The diamond inequality}\label{subsec:sol-diamond}
Recall that $p=-\det(v_0,v_n)$ is the first coordinate of the last rod vector $v_n=(p,q)$. Its absolute value is the area of the
parallelogram spanned by $v_0$ and $v_n$, and by Lemma~\ref{lem:rods-end} it is the order $|\Gamma|$ of the group at infinity. Since
$p\ne0$, the vectors $v_0$ and $v_n$ form a basis of $\R^2$, and we write
\begin{equation}\label{eq:sol-components}
v=\alpha v_0+\beta v_n,\qquad \alpha=-\frac{\det(v,v_n)}p,\qquad \beta=-\frac{\det(v_0,v)}p ,
\end{equation}
for $v\in\R^2$, and $v_j=\alpha_jv_0+\beta_jv_n$ for the rod vectors. Recall that $\varepsilon=(-1)^r$. By Lemma~\ref{lem:sol-degree}, the order $s$ of
the Hankel determinants in Theorem~\ref{thm:soliton} is $r$. The following comparison of components is taken from the proof of \cite[Cor.~3.2(d)]{N4}.

\begin{lemma}\label{lem:sol-frame}
Let $(M,g)$ be a toric ALE instanton with the degree $r$, and let $\mu$ be as in Theorem~\ref{thm:soliton}. Then for $0\le j\le n$,
\begin{equation}\label{eq:sol-normalised}
\frac{H_r(\mu_j)}{H_r(\mu)}=\alpha_j+\varepsilon\beta_j,\qquad \frac{H_{r+1}(\mu_j)}{H_{r+1}(\mu)}=\alpha_j-\varepsilon\beta_j .
\end{equation}
\end{lemma}

\begin{proof}
Let $\mathcal A$ be the linear map that sends $v\in\mathfrak t$, with the components $(x_0,x_1)$ with respect to $(f_0,f_1)$, to
$\mathfrak c^{-1/2}\bigl(x_0/H_r(\mu),x_1/H_{r+1}(\mu)\bigr)$. By \eqref{eq:sol-rodvectors}, $\mathcal Av_j$ is the vector of the two
quotients in \eqref{eq:sol-normalised}. By Lemma~\ref{lem:sol-signs}(i), $\mathcal Av_0=(1,1)$ and $\mathcal Av_n=\varepsilon(1,-1)$.
Since $\mathcal A$ is linear, $\mathcal Av_j=\alpha_j(1,1)+\beta_j\varepsilon(1,-1)$, which is \eqref{eq:sol-normalised}.
\end{proof}

\begin{definition}\label{def:diamond}
For a chain with $p\ne0$ the \emph{open diamond} of its end is the open parallelogram with the vertices $\pm v_0$ and $\pm v_n$,
\begin{equation}\label{eq:sol-diamond}
D(v_0,v_n)=\bigl\{\alpha v_0+\beta v_n:\ |\alpha|+|\beta|<1\bigr\}=\bigl\{v\in\R^2:\ |\det(v_0,v)|+|\det(v,v_n)|<|p|\bigr\}.
\end{equation}
\end{definition}

The two descriptions agree by \eqref{eq:sol-components}. The open diamond is convex and invariant under $v\mapsto-v$. Since
$|a|+|b|=\max(|a+b|,|a-b|)$, it is the set where $|\alpha+\beta|<1$ and $|\alpha-\beta|<1$.

The end rod vectors $v_0$ and $v_n$ are vertices of the diamond, where the sum in \eqref{eq:sol-diamond} equals $|p|$. The diamond
inequality below concerns only the interior rod vectors $v_1,\dots,v_{n-1}$.

\begin{proposition}[the diamond inequality, {\cite[Cor.~3.2(d)]{N4}}]\label{prop:diamond}
Let $(M,g)$ be a toric ALE instanton with $n$ turning points and the degree $r$, and suppose that $1\le r\le n-2$. Then the interior
rod vectors lie in the open diamond of the end:
\begin{equation}\label{eq:sol-diamond-ineq}
|\det(v_0,v_j)|+|\det(v_j,v_n)|<|p|,\qquad 0<j<n .
\end{equation}
\end{proposition}

\begin{proof}
Both $r$ and $r+1$ lie in $\{1,\dots,n-1\}$. So Lemma~\ref{lem:sol-signs}(ii), applied to the orders $r$ and $r+1$, and Lemma~\ref{lem:sol-frame} give
$|\alpha_j+\varepsilon\beta_j|<1$ and $|\alpha_j-\varepsilon\beta_j|<1$ for $0<j<n$, that is, $|\alpha_j|+|\beta_j|<1$.
\end{proof}

In the degrees $r=0$ and $r=n-1$ excluded in Proposition~\ref{prop:diamond} the metric is hyperk\"ahler
(Proposition~\ref{prop:degree0}). For $r=0$ the first quotient in \eqref{eq:sol-normalised} is $H_0(\mu_j)/H_0(\mu)=1$, so
$\alpha_j+\beta_j=1$, and the second one has absolute value at most $1$ by Lemma~\ref{lem:sol-signs}, so $|\alpha_j-\beta_j|\le1$. Hence
$\alpha_j,\beta_j\ge0$, and all rod vectors lie on the segment from $v_0$ to $v_n$, on the boundary of the diamond. For the chains
$(2,\dots,2)$, for instance, $v_j=(j,1-j)=(1-\frac jn)v_0+\frac jnv_n$.

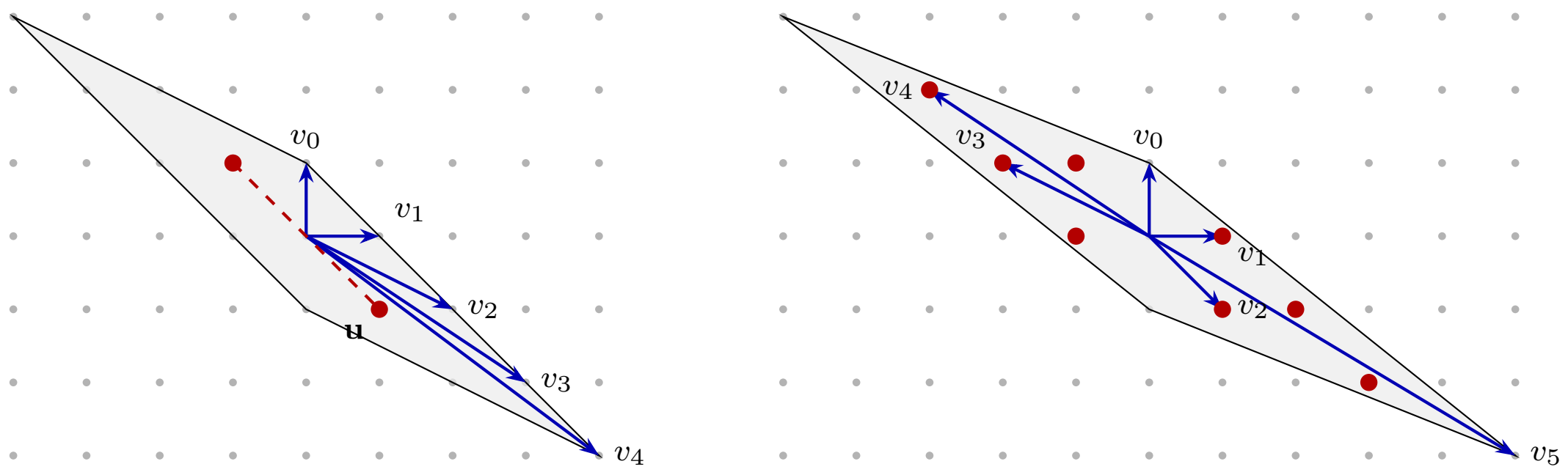
\begin{figure}[!htb]
\centering
\begin{tikzpicture}[x=0.66cm,y=0.66cm,dot/.style={circle,fill=gray!60,inner sep=0pt,minimum size=2pt},
  rod/.style={-{Stealth[length=5pt]},thick,blue!70!black},hub/.style={thick,dashed,red!70!black}]
\begin{scope}
\foreach \x in {-4,...,4}{\foreach \y in {-3,...,3}{\node[dot] at (\x,\y) {};}}
\fill[gray!15,opacity=0.7] (0,1) -- (4,-3) -- (0,-1) -- (-4,3) -- cycle;
\draw[thin] (0,1) -- (4,-3) -- (0,-1) -- (-4,3) -- cycle;
\draw[rod] (0,0) -- (0,1); \draw[rod] (0,0) -- (1,0); \draw[rod] (0,0) -- (2,-1); \draw[rod] (0,0) -- (3,-2); \draw[rod] (0,0) -- (4,-3);
\draw[hub] (-1,1) -- (1,-1);
\foreach \x/\y in {1/-1,-1/1}{\fill[red!70!black] (\x,\y) circle (2.2pt);}
\node[font=\scriptsize,above] at (0,1) {$v_0$};
\node[font=\scriptsize,right] at (4,-3) {$v_4$};
\node[font=\scriptsize,above right] at (1,0) {$v_1$};
\node[font=\scriptsize,right] at (2,-1) {$v_2$};
\node[font=\scriptsize,right] at (3,-2) {$v_3$};
\node[font=\scriptsize,below left] at (1,-1) {$\mathbf u$};
\node[font=\small] at (0,-4.2) {(a) $\kappa=(2,2,2)$};
\end{scope}
\begin{scope}[xshift=7.6cm]
\foreach \x in {-5,...,5}{\foreach \y in {-3,...,3}{\node[dot] at (\x,\y) {};}}
\fill[gray!15,opacity=0.7] (0,1) -- (5,-3) -- (0,-1) -- (-5,3) -- cycle;
\draw[thin] (0,1) -- (5,-3) -- (0,-1) -- (-5,3) -- cycle;
\draw[rod] (0,0) -- (0,1); \draw[rod] (0,0) -- (1,0); \draw[rod] (0,0) -- (1,-1); \draw[rod] (0,0) -- (-2,1); \draw[rod] (0,0) -- (-3,2);
\draw[rod] (0,0) -- (5,-3);
\foreach \x/\y in {1/0,1/-1,2/-1,3/-2,-1/0,-1/1,-2/1,-3/2}{\fill[red!70!black] (\x,\y) circle (2.2pt);}
\node[font=\scriptsize,above] at (0,1) {$v_0$};
\node[font=\scriptsize,right] at (5,-3) {$v_5$};
\node[font=\scriptsize,below right,yshift=1pt] at (1,0) {$v_1$};
\node[font=\scriptsize,right] at (1,-1) {$v_2$};
\node[font=\scriptsize,above left] at (-2,1) {$v_3$};
\node[font=\scriptsize,left] at (-3,2) {$v_4$};
\node[font=\small] at (0,-4.2) {(b) $\kappa=(1,-1,1,-1)$};
\end{scope}
\end{tikzpicture}
\caption{The open diamond $D(v_0,v_n)$ (shaded) with the rod vectors (arrows) and its non-zero lattice points (red dots), for two chains of Table~\ref{tab:chains}. The end rod vectors $v_0$ and $v_n$ are vertices of the diamond. (a) For the
chain $(2,2,2)$ the rod vectors lie on the boundary of the diamond, on the segment from $v_0$ to $v_4$, and the only
non-zero lattice points inside are $\pm\mathbf u=\pm(1,-1)$, joined by a dashed segment, and $\mathbf u$ is not a rod vector. (b) For the chain $(1,-1,1,-1)$ the diamond contains
eight non-zero lattice points on four lines, and all four interior rod vectors lie inside.}\label{fig:diamond}
\end{figure}

\begin{example}\label{ex:diamond}
We compare the diamonds of two ends (Figure~\ref{fig:diamond}). For $v=(x,y)$ we have $\det(v_0,v)=-x$.

The chain $(2,2,2)$ of Table~\ref{tab:chains} has $p=4>0$. Here $\det(v,v_4)=-(3x+4y)$, so the open diamond is the set where
$|x|+|3x+4y|<4$. Its only lattice points are $0$ and $\pm(1,-1)$, all on one line. By Lemma~\ref{lem:lattice} below, no chain with
this end, or with any end $L(m,1)$, and $n\ge3$ has all its interior rod vectors in the open diamond.

The chain $(1,-1,1,-1)$ of Table~\ref{tab:chains} has the end $L(5,2)$ and the degree $r=2$, so that $1\le r\le n-2$. It is invariant
under the change $\kappa\mapsto-(\kappa_{n-1},\dots,\kappa_1)$ that reverses the orientation of $M$ (Section~\ref{subsec:rods-words}).
Here $\det(v,v_5)=-(3x+5y)$, so the open
diamond is the set where $|x|+|3x+5y|<5$. Its non-zero lattice points are $\pm(1,0)$, $\pm(1,-1)$, $\pm(2,-1)$, and $\pm(3,-2)$, on four
different lines. The interior rod vectors $v_1,\dots,v_4$ give the sums $|\det(v_0,v_j)|+|\det(v_j,v_5)|=4,3,3,4<5$, so they lie in
the open diamond, and the chain satisfies \eqref{eq:sol-diamond-ineq}. So the diamond inequality does not exclude a toric ALE
instanton with this chain, and the argument of Section~\ref{sec:proof} says nothing about the end $L(5,2)$. We show this chain because
the chains of toric ALE instantons that are hyperk\"ahler for neither orientation have entries in $\{-1,0,1\}$
\cite[Lemma~A.20]{LiSun2026toric}, \cite[Thm.~2.3]{N4}, and such a chain with $|p|=5$ has at least four entries. Shorter chains of the
end $L(5,2)$ also satisfy \eqref{eq:sol-diamond-ineq} with $1\le r\le n-2$, for instance the chain $(2,-2)$ of Table~\ref{tab:chains}.
\end{example}

\section{Proof of the main theorem}\label{sec:proof}

The proof combines the diamond inequality with an elementary fact about lattice points. If $q\equiv\pm1\pmod{|p|}$, which by
Lemma~\ref{lem:rods-end} means that the group at infinity is conjugate in $\mathrm O(4)$ to a subgroup of $\SU(2)$, then all lattice
points of the open
diamond of the end lie on one line. Two consecutive rod vectors form a lattice basis, so they cannot both lie on this line. By
Proposition~\ref{prop:diamond}, the degree of a toric ALE instanton with such an end is therefore $0$ or $n-1$. In these two
degrees the metric is a multi-Eguchi--Hanson metric (Proposition~\ref{prop:degree0}). For $r=0$ this can also be read off from the
closed formula \eqref{eq:sol-closed} of Theorem~\ref{thm:soliton}, which is then the Gibbons--Hawking form of the metric
\cite[Prop.~11.3]{GerochCF}.

\subsection{The lattice lemma}\label{subsec:pfA-lattice}
We first determine the lattice points of the open diamond of an end with $q\equiv\pm1\pmod{|p|}$, and then show that such a diamond
cannot contain two consecutive rod vectors.

\begin{lemma}\label{lem:diamondpoints}
Let $\kappa$ be a chain with $p\ne0$, and suppose that $q\equiv\sigma\pmod{|p|}$ with $\sigma\in\{1,-1\}$. Then $v_n-\sigma v_0=p\,\mathbf u$
for a primitive vector $\mathbf u\in\Z^2$, and the lattice points of the open diamond $D(v_0,v_n)$ are the multiples $k\mathbf u$ with
$2|k|<|p|$. In particular all of them lie on the line $\R\mathbf u$.
\end{lemma}

\begin{proof}
Since $v_0=(0,1)$ and $v_n=(p,q)$, the vector $v_n-\sigma v_0=(p,q-\sigma)$ is $p$ times $\mathbf u=(1,(q-\sigma)/p)\in\Z^2$, which is
primitive. Let $v=\alpha v_0+\beta v_n$ be a lattice point of the open diamond, so that $|\alpha|+|\beta|<1$. By \eqref{eq:sol-components}
the numbers $l=p\alpha$ and $k=p\beta$ are integers, and $|l|+|k|<|p|$. The second coordinate of $v=(lv_0+kv_n)/p$ is $(l+kq)/p$, which is
an integer, so $l\equiv-kq\equiv-k\sigma\pmod{|p|}$. Since $|l+k\sigma|\le|l|+|k|<|p|$, this forces $l=-k\sigma$. Hence
$v=k(v_n-\sigma v_0)/p=k\mathbf u$, with $2|k|=|l|+|k|<|p|$. Conversely every such $k\mathbf u$ lies in the open diamond.
\end{proof}

For the chains $(2,\dots,2)$ the vector $\mathbf u=(1,-1)$ satisfies $\det(\mathbf u,v_j)=1$ for all their rod vectors (Table~\ref{tab:chains}), so its
slope $-1$ is adjacent in the Farey graph to the slopes of all rod vectors (Remark~\ref{rem:farey} and Figure~\ref{fig:diamond}(a)).
In general, the condition
$q\equiv\sigma\pmod{|p|}$ says that the slope of $\mathbf u$ is adjacent to the slopes of $v_0$ and $v_n$, since $\det(v_0,\mathbf u)=-1$
and $\det(\mathbf u,v_n)=\sigma$. Lemma~\ref{lem:diamondpoints} says that in the open diamond only the line $\R\mathbf u$ is available, which by the
opening of this section is too little for two consecutive rod vectors.

\begin{lemma}[lattice lemma]\label{lem:lattice}
Let $\kappa$ be a chain with $n\ge3$ and $p\ne0$ whose interior rod vectors $v_1,\dots,v_{n-1}$ lie in the open diamond $D(v_0,v_n)$.
Then $q\not\equiv\pm1\pmod{|p|}$. In particular $|p|\ge5$ and $|p|\ne6$.
\end{lemma}

\begin{proof}
Suppose that $q\equiv\pm1\pmod{|p|}$. By Lemma~\ref{lem:diamondpoints}, $v_1$ and $v_2$ are multiples of the same vector $\mathbf u$, so
$\det(v_1,v_2)=0$, which contradicts \eqref{eq:rods-det}. The last statement holds since every unit modulo $m\in\{1,2,3,4,6\}$ is
congruent to $\pm1$ (Section~\ref{subsec:rods-end}).
\end{proof}

\subsection{The degrees \texorpdfstring{$0$}{0} and \texorpdfstring{$n-1$}{n-1}}\label{subsec:pfA-explicit}
In the two extreme degrees the metric is hyperk\"ahler for one of the orientations, by a theorem of Li and Sun, and the toric
hyperk\"ahler ALE instantons are known. The following proposition is \cite[Cor.~3.2(c), Rem.~3.6]{N4}.

\begin{proposition}\label{prop:degree0}
Let $(M,g)$ be a toric ALE instanton with $n$ turning points and the degree $r=0$. Then its chain is $(2,\dots,2)$, $p=n$, and
$(M,g)$ is a multi-Eguchi--Hanson metric with $n$ collinear centres that is hyperk\"ahler for the orientation opposite to that of $M$.
If instead $r=n-1$, then the chain is $(-2,\dots,-2)$, $p=(-1)^{n+1}n$, and $(M,g)$ is a multi-Eguchi--Hanson metric with $n$ collinear centres that is
hyperk\"ahler for the orientation of $M$.
\end{proposition}

\begin{proof}
Let $r=n-1$. By Remark~\ref{rem:LiSun-degree}, the degree of Li and Sun of $(M,g)$ is $\rb$, and $\rb=n-1-r=0$ by
Lemma~\ref{lem:rods-degree}. By the theorem of Li and Sun
\cite[Thm.~7.4]{LiSun2025}, in the form of \cite[Thm.~2.2]{N4}, $g$ is therefore hyperk\"ahler for the orientation of $M$, and it is a
Gibbons--Hawking metric with collinear centres. The centres are the $n$ turning points \cite[\S11.2]{GerochCF}, and since the end is
ALE, $g$ is a multi-Eguchi--Hanson metric \cite[Rem.~3.6]{N4}. That remark is stated for $n\ge2$, and for $n=1$ a Gibbons--Hawking
metric with one centre and an ALE end is flat, the multi-Eguchi--Hanson metric with one centre. The chain is $(-2,\dots,-2)$ by
\cite[Cor.~3.2(c)]{N4}, since the degree index of \cite{N4} is $s=r=n-1$ (Lemma~\ref{lem:sol-degree}), and $p=(-1)^{n+1}n$ by Table~\ref{tab:chains}. If $r=0$, then $M$ with the opposite orientation has the degree $n-1$, since reversing the orientation exchanges $r$ and
$\rb$ (Section~\ref{subsec:rods-words}), and the case just treated applies to it. So $g$ is a multi-Eguchi--Hanson metric that is hyperk\"ahler for the orientation opposite to that of
$M$, the chain of $M$ is $-(-2,\dots,-2)=(2,\dots,2)$, and $p=n$ by Table~\ref{tab:chains}.
\end{proof}

\subsection{The main theorem}\label{subsec:pfA-main}
We now combine the diamond inequality, the lattice lemma, and the description of the extreme degrees.

\begin{theorem}\label{thm:main}
Let $(M,g)$ be a toric ALE instanton whose group at infinity $\Gamma$ is conjugate in $\mathrm O(4)$ to a subgroup of $\SU(2)$. Then
the action has exactly $|\Gamma|$ fixed points. For one of the orientations of $M$, $(M,g)$ is a multi-Eguchi--Hanson metric with $|\Gamma|$ collinear centres that is hyperk\"ahler for this
orientation, and its chain is $(-2,\dots,-2)$. For $|\Gamma|=1$ it is the flat $\R^4$, and for $|\Gamma|=2$ it is the Eguchi--Hanson
metric.
\end{theorem}

\begin{proof}
By Lemma~\ref{lem:basic}, $M$ is simply connected, the stabilisers are connected, and the action has $n\ge1$ fixed points. Fix
orientations of $M$ and $\mathfrak t$, and let $\kappa$ be the chain, $(p,q)=v_n$, and $r$ the degree. By Lemma~\ref{lem:rods-end},
$|p|=|\Gamma|$, and the hypothesis says that $q\equiv\pm1\pmod{|p|}$ (Section~\ref{subsec:rods-end}). Suppose that
$1\le r\le n-2$. Then $n\ge3$, and by Proposition~\ref{prop:diamond} the interior rod vectors lie in the open diamond $D(v_0,v_n)$.
By Lemma~\ref{lem:lattice}, $q\not\equiv\pm1\pmod{|p|}$, which contradicts the hypothesis. Hence $r=0$ or $r=n-1$, and
Proposition~\ref{prop:degree0} shows that $(M,g)$ is a multi-Eguchi--Hanson metric with $n$ collinear centres, with the chain
$(-2,\dots,-2)$ in its hyperk\"ahler orientation. In both cases $n=|p|=|\Gamma|$. A multi-Eguchi--Hanson metric with one centre is
flat, and with two centres it is the Eguchi--Hanson metric.
\end{proof}

The hypothesis on $\Gamma$ enters only through the lattice lemma.

\begin{corollary}\label{cor:smallgroups}
Let $(M,g)$ be a toric ALE instanton whose group at infinity has order $|\Gamma|\in\{1,2,3,4,6\}$. Then $(M,g)$ is a
multi-Eguchi--Hanson metric with $|\Gamma|$ collinear centres.
\end{corollary}

\begin{proof}
By Lemma~\ref{lem:rods-end}, $|\Gamma|=|p|$, and $q$ is a unit modulo $|p|$. For $|p|\in\{1,2,3,4,6\}$ every unit modulo $|p|$ is congruent to
$\pm1$ (Section~\ref{subsec:rods-end}), so $q\equiv\pm1$, and Theorem~\ref{thm:main} applies.
\end{proof}

Our argument does not cover $|\Gamma|=5$. In that case, the end $L(5,2)$ is possible a priori, and for this end the diamond inequality
alone does not suffice (Example~\ref{ex:diamond}). Computer calculations by the authors, which use a reduction to finitely many chains
and interval arithmetic and are not part of this paper, indicate that the conclusion of Corollary~\ref{cor:smallgroups} holds for
$|\Gamma|=5$ as well.
\section*{Computation and the use of an AI assistant}
The project was carried out with the help of an AI assistant, Claude (Anthropic), used through the Claude Code environment. The
assistant took part in developing the arguments and in writing the text. We also used it to review the literature and to read the
drafts critically. Computer calculations were used only for the remark on the order $|\Gamma|=5$ at the end of Section~\ref{sec:proof}. The proofs in this paper do not use them. The authors take responsibility for the content.

\bibliographystyle{amsplain}
\bibliography{refs}

\end{document}